\documentstyle{amsppt}
\magnification=1200
\NoBlackBoxes
\def\ss{\vskip.10in}
\def\ls{\vskip.25in}
\def\P{\Bbb P}
\def\C{\Bbb C}

\def\cc{\Cal C}
\def\Z{\Bbb Z}

\centerline{\bf Bend, Break and Count}
\centerline{Preliminary Version: 27.2.97}
\vskip.25in
\centerline{Z. Ran}
\centerline{Department of 
Mathematics} 
\centerline{Durham University, UK and University of California, Riverside}
\centerline{Riverside, CA  92521 USA}
\centerline{ziv\@ucrmath.ucr.edu}
\vskip.5in
The purpose of this paper is establish and apply an enumerative formula
or `method' dealing with a family $\cc = \{\bar{C}_y : y \in Y\}$ of rational 
curves on a variety $S$, e.g. a rational surface.  Significantly, the 
family $\cc$ is not assumed to be the family of `all' rational curves of given
homology class: rather, we require only that it be sufficiently large
$(n = \dim Y \geq 3)$ and well-behaved as regards deformation theory
and codimension-1 degenerations.  The formula computes the `degree' 
$d(\cc)$, i.e. the number of curves $\bar{C}_y$ through $n$ general
points of $S$, in terms of analogous degree-type numbers attached to the 
(codimension-1) boundary component's $Z$ of $Y$, which parametrize the 
reducible curves (there are several such numbers depending on how the 
components of the reducible curve are `weighted').

Some comments are in order on connections with quantum cohomology.  While
the author denies any first-hand knowledge or understanding of the latter, its algebro-geometric aspect has been represented as essentially equivalent
to certain recursive formulae for counting rational curves, which are 
contained in the associativity formula for quantum multiplication.
It has seemed to the author that these recursions, at least, are 
largely a matter of taking advantage of the `slack' in the problem,
i.e. the large number of deformation parameters for rational curves 
(on `convex' varieties).  This viewpoint suggests a connection
with Mori's bend-and-break technique, a small part of which is the 
observation that once a rational curve `bends'  sufficiently (on a 
surface, this means moving in a 3-parameter family) it will `break',
i.e. admit a reducible limit.  Our general formula (Theorem 1, Sect. 1)
is merely a quantitative version of this idea.  As already indicated, it 
applies to any given(good enough) family and accordingly does not 
rely on existence of (compact) moduli spaces for (reparametrization classes
of) stable maps as in [FP].  The proof is a completely elementary 
argument involving (multi-) sections and fibre components on a birationally
ruled surface.

Now for better or worse, the effect of Theorem 1 is to shift the difficulty
elsewhere, namely to `enumerating' the boundary components $Z$, which in 
principle is a lower-degree problem, but not necessarily well-behaved.
The simplest $Z$ are of `product type' and parametrize a pair of independently
varying curves plus a point of their intersection: these are unproblematic.
However, there are others, such as those parametrizing a pair of 
mutually tangent curves, and worse: a variety of examples is given in 
Sect. 2.

In Sect. 3 we consider the problem of enumerating plane curves of given 
degree $d$ and given moduli, i.e., birational to a fixed smooth curve 
$C$ of genus $g$ (the analogous problem with fixed $d,g$ and unrestricted
moduli, a.k.a. the degree of the Severi-variety having long been settled
[R]).  If $N(d,g\rangle$ denotes the number of such curves 
through $3d-2g + 2$ general points (or $3d-1$ if $g=1$) then Pandharipande
[P] has shown
$$
N(d,1\rangle = \frac{(d-1)(d_2)}{2} N(d,0) = \frac{(d-1)(d-2)}{2} N_d.
$$
Here we give a recursive procedure for computing $N(d,2\rangle$.  While 
it is fairly clear what sub-problems would have to be solved to extend 
this to higher genus, it is unclear whether those can be solved, especially 
ones 
involving high-order contact between rational curves.  However, see [R$'$]
for a different approach to this problem. 

\subheading{Acknowledgment}  This work was begun at the Mittag-Leffler
Institute while the author was a participant in its program on 
Enumerative Geometry, led by W. Fulton.  It is a pleasure to thank the 
Institute's staff, the other participants, and most especially Bill Fulton   
for their contributions, which have made this possible.  

\subheading{1. General formula}

Suppose given a flat family  $\{\bar{C}_y : y \in Y\}$ of curves on a 
smooth projective surface $S$, parametrized by an irreducible $n$-dimensional
projective variety $Y$, such that a general $\bar{c}_y$ is a irreducible 
rational curve.   Thus we have a diagram
$$
\aligned
&\cc \subset Y \times S {\overset f=p_2 \to\rightarrow} S \\
\bar{\pi}&\downarrow\\
&Y
\endaligned
\tag1.1
$$
where $\bar{\pi}$ is flat with fibres $\pi^{-1} (y) = \bar{C}_y$.  We
assume $Y$ is normal.  In what follows, only geometry is codimension
$\leq 1$ on $Y$ will play any role (so we could actually assume $Y$ 
nonsingular without either losing generality or gaining convenience).  Let 
$n: \cc \to \bar{\cc}$ be the normalization and $\pi = \bar{\pi}$ on the resulting
flat family with fibres $C_y = \pi^{-1} (y)$ mostly isomorphic to 
$\P^1$.  We let $\partial Y  \subset Y$, the `boundary', denote the 
discriminant locus of $\pi$.  We now introduce a strong `good behavior'
condition on our family which, while not absolutely essential, makes
for a simpler enumeration formula and is satisfied in applications.

\noindent
\roster
\item{$(*)$}  $Y$ is normal, hence $\cc$ is smooth in codimension 2 except at 
singular points of fibres; for a generic point $z$ of any 
$(n-1)$-dimensional component $Z$ of $\partial Y$, the fibre $C_z$ has 
just two components $C_{1,z}, C_{2,z}$, both $\P^1$'s, and $C_{1,Z} \cap
C_{2,Z} = \{p\}$ is an $A_{\ell-1} \times \C^{n-1}$ singularity on $e$
for some $\ell = \ell (Z)$.  
\endroster

Now the {\it degree} (classically, the {\it grade}) of our family, which 
we abusively denote by $d(Y)$, is defined to be the number of $y\in Y$ 
such that $\bar{C}_y$ contains $n = \dim Y$ general points $s_1, \ldots,
s_n \in S$.  To define analogues of degree for boundary components, it
will be convenient to count each one twice by introducing a `making'.
Thus define a {\it marked} boundary component $Z$ to consist of a boundary
component $Z'$ together with an ordering of the two corresponding
components $\bar{C}_{1,Z}, \bar{C}_{2,Z}$; if the components are 
monodromy-interchangeable we still consider both orderings.  To this 
data we associate the following degrees
$$
\align
d^{1,1} (Z) &= \# \{Z : \bar{C}_z \ni s_1, \ldots, s_{n-1}, 
s, \in \bar{C}_{1,z} s_2 \in \bar{C}_{2,z} \} \\
d^{0,2} (Z) &= \# \{ Z : \bar{C}_z \ni s_1, \ldots, s_{n-1} , 
s_1, s_2 \in \bar{C}_{2,z} . 
\tag1.2
\endalign
$$

The simplest-though by no means all -- boundary loci are those where 
the two components vary independently.  More precisely, let us say that a 
boundary locus, sum of boundary components, is of {\it product type}
if it can be naturally identified with a locus $\{(\bar{C}_1, \bar{C}_2,
p \in \bar{C}_1 \cap \bar{C}_2)\}$, where $\bar{C}_i \in 
\{ \bar{C}_i\}$ are independtly generic in their respective
(irreducible) families of dimension $n_i, n_1 + n_2 = n-1, \bar{C}_1$
and $\bar{C}_2$ meet transversely and $p \in \bar{C}_1 \cap \bar{C}_2$
can be specified arbitrarily.  For $Z$ of product type, clearly
$$
\align
d^{1,1} (Z) &= \pmatrix n-1\\ n_1 - 1\endpmatrix d(\{\bar{C}_1\}) 
d(\{\bar{C}_2\}) \bar{C}_1 . \bar{C}_2\\
d^{0,2} &= \pmatrix n-1\\ n_1\endpmatrix d(\{C_1\}) d(\bar{C}_2\})
\bar{C}_1 . \bar{C}_2.
\tag1.3
\endalign
$$

\proclaim{Theorem 1}  In the above situation, suppose moreover $n\geq 3$.
Then for any line bundle $L$ on $S$, we have
$$
L^2 d(Y) = \sum \ell (Z) [ (C_1 L) (C_2 L) d^{1,1} (Z) - 
(C_1 L)^2 d^{0,2} (Z) ] , 
\tag1.4
$$
the sum being over all marked boundary component $Z$ with corresponding
curves $C_1, C_2$
\endproclaim

\demo{Proof}  First we might as well cut $Y$ down to a 3-fold by 
imposing $s_{4}, \ldots, s_n$ and henceforth assume $n=3$.  Then cut 
$Y$ down to a (smooth) curve $B$ by imposing $s_1, s_2$ and 
minimally resolve $\cc_Y \times B$, thus obtaining a diagram
$$
\aligned
&X\quad {\overset f\to \rightarrow}\quad S\\
\pi &\downarrow\\
&B
\endaligned
\tag1.5
$$
with $X$ a smooth surface, $\pi$ a blown-up $\P^1$-bundle with 
sections $S_i$ corresponding to $s_i, i = 1,2$, and reducible fibres
of the form
$$
C^i = C^i_1 + E^i_1 + \cdots + E^i_{\ell_i - 1} + C_2^i
\tag1.6
$$
with $f(E_1^i + \cdots + E^i_{\ell-1}) = p^i$, a point on $S$, and 
$\ell_i = \ell (Z_i)$ where $Z_i$ is the boundary component whence 
$C^i$ comes.  For future use, we note here that knowing the 
structure of the reducible fibres $C^i$ is equivalent to knowing the  
singularity type of $\cc$ along $Z_i$, and in any event is the only thing 
we need to know this singularity for.  It is obvious-but 
nevertheless crucial - that the Neron-Severi group of $X$ is generated
by any section plus fibre components.  Note that $(E_j^i)^2 = - 2, 
(C_j^i)^2 = - 1$.  It is easy to see from this firstly that 
$$
S_1 - S_2 \sim \sum\Sb s_1 \in C_1^i\\ s_2 \in C_2^i\endSb \left(
\sum_{j=1}^{\ell_i - 1} (jE_j^i + \ell_i C_2^i) + 
\sum\Sb s_2 \in C_1^i\\ s_1 \in C_2^i\endSb \left( \sum_{j=1}^{\ell_1 -1}
(\ell_i - j) E_j + \ell_i C_1^i\right)\right) - m F, 
\tag1.7
$$
$F$ = fibre.  Taking the dot product with $S_1$ yields $-m = S_1^2$.

As $s_1$ and $s_2$ are interchangeable by a suitable monodromy
tranformation and dot products are preserved, we also have $-m = S_2^2$, so 
taking the dot product of (1.7) with $S_2$ yields
$$
m = \sum \Sb s_1 \in C_1^i\\ s_2 \in C_2^i \endSb \ell_i = \sum_Z \ell(Z) 
d^{1,1} (Z). 
\tag1.8
$$
Now set $d = L. \bar{C}_y$ and note as before that $dS_1 - f^* L$ is a 
linear combination of fibre components, hence one can easily check that 
$$
dS_1 - f^*L \sim \sum\Sb s_1 \in C_1^i \\ s_2 \in C_1^i v C_2^i \endSb 
(C_2^i . L) (\sum_{j=1}^{\ell_i - 1} jE_j^i + \ell_i C_2^i) + 
\sum\Sb s_1 \in C_2^i \\ s_2 \in C_1^i \cup C_2^i \endSb 
(C_1^i . L) (\sum (\ell_i - j) E_j^i + \ell_i C_1^i ) 
- aF 
\tag1.9
$$
Then taking dot product with $S_1$ yields $a = dm$.  Then squaring both 
sides we get
$$
\align
L^2 d(Y) &= (f^* L)^2 = d^2m - \sum\Sb s_1 \in C_1^i \\ s_2 \in C_1^i v
C_2^i \endSb \ell_i (C_2^i.L)^2 - 
\sum\Sb s_1 \in C_2^i \\ s_2 \in C_1^i \cup C_2^i \endSb \ell_i
(C_1^i . L)^2 \\
&= \sum\Sb s_1 \in C_1^i \\ s_2 \in C_2^i\endSb (C_1^i . L + C_2^i.L)^2 
\ell_i -  \sum\Sb s_1 \in C_1^i \\ s_2 \in C_1^i\cup 
C_2^i \endSb \ell_i (C_2^i.L)^2 - 
\sum\Sb s_1 \in C_2^i \\ s_2 \in C_1^i v C_2^i \endSb \ell_i
(C_1^i . L)^2 \\
&= 2 \sum\Sb s_1 \in C_1^i \\ s_2 \in C_2^i \endSb (C_1^i . L) (C_2^i . L) 
\ell_i - \sum_{s_1, s_2 \in C_1^i} \ell_i (C_2^i . )^2 - 
\sum_{s_1, s_2 \in C_2^i} \ell_i (C_1^i . L)^2\\
&= \sum d^{1,1} (Z) \ell (Z) (C_1 . L)(C_2.L) - \sum d^{0,2} (Z)
\ell (Z) (C_{1,z} . L)^2
\endalign
$$
\qed
\enddemo

\subheading{2.  Various Examples}
\ss
{\it (a)  Del Pezzo surfaces: all curves}
\ss
We begin by recalling some facts.  Let $S$ be a Del Pezzo surface, with 
(ample) anticanonical bundle $-K$ and Neron-Severi group $N$.  Let us 
call a class $C\in N$ {\it good} if $-K C \geq 0$ and 
$C^2 \geq 0$.  Then 

(i)  if $C$ is the class of an integral curve $\bar{C}$ then either
$\bar{C}$ is a line $(-K\bar{C} = 1, \bar{C}^2 = - 1$), or $C$ is good;

(ii) If $\bar{C}$ is an irreducible rational curve then as such $\bar{C}$
has unobstructed defomations of dimension $-K . \bar{C} - 1$
which are generally transverse to given subvarieties.

It follows easily that if $\bar{C} \subset S$ is a good rational curve and 
$Y$ is the normalization of the locus $\{\bar{C}\}$ of rational curves 
in the linear system $|\bar{C}|$, then the hypotheses of Theorem 1 are 
satisfied provided $n = - K . \bar{C} - 1 \geq 3$.  The boundary loci 
$Z$ all have $\ell (Z) = 1$, are of product type and correspond to 
(ordered) expressions
$$
c = [\bar{C}] = c_1 + c_2
\tag2.a.1
$$
with $c_1, c_2$ representable by (irreducible) rational curves.  Explicitly,
taking $L = - K$ leads to the following, in which we denote $d(Y)$ by $N_c$:
$$
\align
K^2 N_c & = \sum_{c = c_1 + c_2} N_{c_1} N_{c_2} \left[
(-K. c_1) (-K. c_2) (c_1 . c_2) \pmatrix -K.c - 4\\ -K.c_1 -2 \endpmatrix\right.\\
 &\left.- (-K.c_1)^2 c_1 . c_2 \pmatrix - K.c -4\\ - K.c_1 - 1\endpmatrix \right]
\tag2.a.2
\endalign
$$
Now at least if $S$ has anticanonical degree $\geq 3$, then
any good class $C$ is representable 
by a rational curve
 (see Appendix).  Whenever this is so, the problem of effectively 
computing the RHS of (2.a.2) becomes a purely combinatorial matter.
For $S = \P^2$, (2.a.2) reduces to the `associativity relation' of 
Kontsevich et al., cf. [FP].  
\ss
{\it (b)  Higher dimensions}
\ss
There are potentially several ways to meaningfully extend Theorem 1 to the 
case of a higher-dimensional ambient variety.  Without getting systematically
involved in this matter here, we shall merely indicate a realtively obvious
such extension, obtained by simply replacing point conditions by incidence
with respect to codimension -2 linear spaces.  Let $\{\bar{C}_y : y \in Y\}$
etc. be as in \S1, except that $m = \dim \ S$ is no longer assumed 
to equal 2, while the line bundle $L$ is assumed very ample.  We may then 
define 
$$
d_L (Y) = \# \{ y : \bar{C}_y \cap L_1^i \cap L_2^i \neq \emptyset \quad
i = 1, \ldots, n \}, 
L_j^i \in |L| \qquad \text{general} 
$$
and likewise for the $d_L^{i,j} (Z)$.  The same arguments apply, notwithstanding
that the analogues of $S_1$ and $S_2$ are now only multisections:
the essential point is that, still,  $S_1^2 = S_2^2$.  The following
formula then obtains
$$
d_L (Y) = \sum \ell (Z) [ (C_1 L) (C_2L) d_L^{1,1} (Z) - (C_1.L)^2
d_L^{0,2} (Z) ] 
\tag2.b.1
$$

\ss
{\it (c)  The plane:  some codimension-1 counts}
\ss
Here we give counts associated with some codimension-2 loci in the family
of all rational plane curves, beginning with the number $C_d$ of rational
curves of degree $d$ through $3d-2$ points having a node in a fixed line $M$.
The marked boundary components $(Z, C_1, C_2)$ are easily 
determined and some in two types, depending on whether $C_1$ and $C_2$ have 
a common point on $M$, or whether $C_1$ or $C_2$ has node on $M$.  For the 
first type we have, e.g. 
$$
d^{1,1} (Z) = \left[\pmatrix 3d-5\\ 3d_1 - 2 \endpmatrix d_1 + \pmatrix
3d-5 \\ 3d_2 -2\endpmatrix d_2\right] (d_1 d_2 -1) N_{d_1} N_{d_2} , 
\tag2.c.1
$$
while the second is of product type.  Applying Theorem 1, we obtain a 
recursion:
$$
\align
C_d &= \sum \left[\left(\pmatrix 3d-5\\3d_1 -2\endpmatrix {d_2} - 
\pmatrix 3d-5 \\ 3d_1 - 1\endpmatrix {d_1} \right) d_1^2 \right.\\
&\left.+ \left(\pmatrix 3d-5 \\ 3d_2 -2\endpmatrix d_1 d_2 - \pmatrix 3d-5 \\ 
3d_2 -3
\endpmatrix {d_1^2} \right) d_2 \right] (d_1 d_2 - 1) N_{d_1} N_{d_2} \\
&+ \sum \left[ \left(\pmatrix 3d-5\\ 3d_1 - 3\endpmatrix + \pmatrix
3d-5 \\3d_2 - 3 \endpmatrix \right) d_1 d_2 \right.\\
&\left.- \pmatrix 3d^5 \\ 3d_1 - 2
\endpmatrix d_1^2 - \pmatrix 3d -5\\ 3d_2 - 2\endpmatrix  d_2^2\right]
d_1 d_2 C_{d_1} N_{d_2}
\tag2.c.2
\endalign
$$

Next we consider some numbers which may be derived from $C_d$ by elementary
means.  First let $B_d$ be the number of rational curves of degree $d$
through $3d-2$ general points which are properly (i.e. at a smooth point)
tangent to a fixed line $M$.  Then $B_d$ is related to $C_d$ by the formula
$$
B_d + 2C_d = 2 (d-1) N_d . 
\tag2.c.3
$$
This comes about by considering the (rational) `restriction' map
$$
\align
r: V_{d,0} \subset  \P^n & \cdots \to  \P^d\\
\bar{C} &\mapsto  \bar{C} \cap M 
\endalign
$$
$$
V_{d,0} = \ \ \text{variety of degree  $-d$ rational curves,}
$$
which pulls back 
the discriminant hypersurface (of degree $2(d-1)$) to the sum of the 
properly tangent locus (with multiplicity 1) plus the `node on $M$'
locus (with multiplicity 2).

A natural generalization of $B_d$ is the number $B_{d,e,g}$ of 
curves (rational degree -$d$, through $3d-2$ genral points) properly
tangent to a given curve $E$ of dgree $e$ and geometric genus $g$.
To compute this we return to the situation of (1.5) where now $f: X \to 
\P^2$ has degree $N_d$, and note that $B_{d,e,g}$ coincides with the number 
of proper (smooth) ramification points of $\pi |_{f^{-1} (E)}$, and that the 
singularities of $E$ reduce the geometric genus of $f^{-1} (E)$ by 
$N_d ((e-1) (e-2)/2 \ \ -g)$.  Hence by the adjunction formula
$$
\align
B_{d,e,g} &= e f^* L (e f^* L + K_X - K_B) - 2N_d ((e-1)(e-2)/2 \ \ - g) \\
&= e(e-1) N_d + e f^*L(f^* L+ K_X - K_B) - (e-1)(e-2) N_d + 2g N_d \\
&=2(e-1) N_d + eB_d + 2gN_d
\tag2.c.4
\endalign
$$
where $L=$ line.  In particular for $g=0$ we obtain
$$
B_{d,e,0} = 2(e-1) N_d + eB_d . 
\tag2.c.5
$$

\ss
{\it (d)  some cross-ratio counts}
\ss
Here we consider some (still codimension-1) counts involving a `marked'
curve $\bar{C}$, i.e. the normalization $C = \P^1$ carries a fixed-up 
to isomorphism-quadruple, or more precisely ordered pair of unordered
pairs $(\{p_1, q_1\}\{p_2, q_2\})$.  Note to begin with that 
$$
Aut (\P^1, \{p_1, q_1\} \{p_2, q_2\}) = \Z_2
$$
inducing the permutation $(p_1, q_1) (p_2 q_2)$ (if, e.g. $(p_1 q_1) = 
(0,\infty)$ this is given by $\frac{p_2 q_2}{z}$).  Now fix integral nodal
plane curves $E_1, F_1, E_2, F_2$ in general position of respective degrees
$e_i, f_j$ and `define'.
$$
N (d < (e_1), (e_2)) = \# \{f: (\P^1 , \{p_1, q_1\} , \{p_2, q_2\} ) 
\to (\P^2 , E_1, E_2 ) \} ;
$$
in this and below it is understood that $f(\P^1)$ is to be a rational 
curve of degree $d$ through $3d-2$ general points of $\P^2$, and $f$
is taken up to automorphism.

Similarly let 
$$
N (d < e_1, f_1 , e_2) = \# \{ f: ( \P^1 , p_1, q_1, \{p_2, q_2\})
\to ( \P^2 , E_1, F_1, E_2 ) \}
$$
and likewise $N(d < e_1, f_1, e_2, f_2)$.  It is, in fact, easy to see
by specialization and induction that these numbers depend only on  the 
degrees of the curves in question:  for instance
$$
\multline
N( d < (e_1 + f_1) , (e_2)) =  N(d < (e_1), (e_2)) + N(d  < (f_1) (e_2)) \\
+ N ( d < e_1, f_1, (e_2)) + N (d < f_1, e_1, (e_2)) , 
\endmultline
$$
etc.; also these numbers possess an evident symmetry, e.g.
$$
N ( d < e_1, f_1, (e_2)) = N (d < f_1, e_1, (e_2)).
$$
It follows formally that 
$$
\align
N ( d < (e_1), (e_2)) &= e_1 e_2 N (d < (1), (1)) + (e_1 - 1) (e_2 - 1) 
(e_1 + e_2) N (d < 1,1, (1))\\
&+ (e_1 -1) (e_2 -1) e_1 e_2 N(d < 1, 1, 1, 1), \\
N(d < e_1, f_1, (e_2)) & = e_1 f_1 e_2 (N ( d < 1,1,(1)) + (e_2 -1) N (d < 
1,1,1,1))\\
N (d < e_1, f_2, e_2, f_2) &= e_1 f_1 e_2 f_2 N (d < 1,1,1,1) 
\tag2.d.1
\endalign
$$

It thus suffices to compute the basic numbers 
$$
N (d < (1), (1)) , \quad N ( d < 1,1, (1)), N (d < 1,1,1,1) ,
$$
for which we may set up a recursion in $d$ based on Theorem 1.
Consider, e.g., the case of $N(d < (1), (1))$.  The curve
$C = \P^1 \to \bar{C}$ of degree $d$ will be marked with $\{p_1, q_1\} , 
\{ p_2, q_2\}$, and the boundary components may be determined by an 
easy dimension count, e.g. based on the deformation theory of the moduli
spaces of stable maps [FP] (through this is not essential).  They 
correspond to pairs $(C_1, C_2)$ where $C_1$ is of degree $d$ marked with 
$\{p_1^0, q_1\}, \{p_2, q_2\}$ and isomorphic as such to $C_1, p_1^0 = 
C_1 \cap C_2$ and $C_2$ is of degree $d_2$ and marked with a `new'
$p_1$ playing the role of the old.  It is clear that for such a 
component $Z$ we have
$$
\align
d^{1,1} (Z) & = N (d_1 < 1, d_2, (1)) d_2 
\pmatrix 3d-5 \\ 3d_2 - 2 \endpmatrix + N_{d_1} N_{d_2} d_1^2
(d_1 - 1) \pmatrix 3d-5 \\ 3d_1 - 2 \endpmatrix \\
&= N(d_1 < 1,1, (1)) d_2^2 \pmatrix 3d-5\\ 3d_2 - 2\endpmatrix + 
N_{d_1} N_{d_2} d_1^2 (d_1 -1) \pmatrix 3d-5 \\ 3d_1 - 2 \endpmatrix.
\tag2.d.2
\endalign
$$
Here the first summand comes from choosing $\bar{C}_2$ through 
$1 + 3d_2 - 2$ points, then a point on $\bar{C}_2 \cap E_1$, then 
$\bar{C}_1$; the second from choosing $\bar{C}_1$ through $1 + 3d_1 - 2$
points, then as ordered pair on $\bar{C}_1 \cap E_2$, then a point 
$\bar{q}_1 \in \bar{C}_1 \cap E$ - which in turn determines 
$p_1^0$ via cross-ratio -- then finally $\bar{C}_2$.  Applying Theorem~1
(recall that each $Z$ has two markings so must be counted twice), we 
conclude:
$$
\align
N(d< (1),(1)) &= \sum \left[ d_1d_2^3 \pmatrix 3d-5\\ 3d_2 - 2 \endpmatrix
 - d_1^2 d_2^2 \pmatrix 3d-5 \\ 3d_1 - 1 \endpmatrix \right] N (d_1 < 
1,1, (1))\\
& \left[ d_1^3 d_2 \pmatrix 3d-5\\ 3d_1 -2\endpmatrix - d_1^4
\pmatrix 3d-5 \\ 3d_2 - 1 \endpmatrix \right] N (d_2 < 1,1, (1))\\
&+ \left[ d_1^3 d_2 (d_1 - 1) \pmatrix 3d-5 \\ 3d_1 -2\endpmatrix - 
d_1^4 (d_1 - 1) \pmatrix 3d-5 \\ 3d_1 - 1 \endpmatrix \right.\\
&\left. + d_1 d_2^3 (d_2 -1) \pmatrix 3d-5\\ 3d_2 - 2 \endpmatrix - 
d_1^2 d_2^2 (d_2-1) \pmatrix 3d-5 \\ 3d_2 -1\endpmatrix \right]
N_{d_1} N_{d_2}
\tag2.d.3
\endalign
$$
Recursions for $N(d< 1,1, (1))$ (involving $N(d < 1, 1, (1)), N (d < 
1,1,1,1)$ and $N_d$) and for $N( d < 1,1,1,1)$ (involving
$N(d < 1,1,1,1)$ and $N_d$ may be obtained similarly.
\ss
For subsequent applications we require a `dual' cross-ratio count when 
the marked curve $\bar{C}_1$ and $\bar{E}_1$ are fixed, as is the marking's
cross-ratio, while $\bar{E}_2$ is allowed to vary (of course as a rational
curve of given degree $e_2$ and through $3e_2 -2$ general points).
Thus define numbers
$$
N(d,(e_1) > (e_2)) = \# \{f: (\P^1, \{p_1, q_1\}, \{p_2, q_2\}) \to 
(\bar{C} , \bar{E}_1, \bar{E}_2 ) \}
$$
with the usual provisi, where $\bar{C}$ is fixed rational degree-$d$, $E_1$
fixed integral nodal of degree $e_1, \bar{E}_2$ rational degree-$e_2$
through $3e_2 -2$ points, and the source of $f$ is fixed up to isomorphism.  We 
analogously define numbers $N(d, e_1, f_1 > (e_2)) , N(d, (e_1), f_2 > e_2),
N(d, e_1, f_1, f_2 > e_2 )$.  As before these behave simply with 
respect to the fixed unmarked curves, e.g.,
$$
N(d, (e_1) > (e_2)) = e_1 N (d, (1) > (e_2)) + 
e_1 (e_1 -1) N (d, 1,1 > (e_2)) . 
\tag2.d.4
$$

Now to compute, e.g., $N(d, 1,1, (e_2))$, the method of Theorem 1 yields
a recursion in $e_2$.  The boundary components $Z = \{E_{2,1} \cup
E_{2,2}\}$ are easily determined and fall into two types depending on 
whether $p_2, q_2 \in E_{2,1}$, say (type (1)) or $p_2 \in E_{2,1},
q_2 \in E_{2,2}$ (type (2)).  The type (1) components are easily 
enumerated in terms of $N(d, 1, 1 > (e_{2,1}))$ and the type (2)'s
in terms of $N(d, 1,1, e_{2,1} > e_{2,2}) = e_{2,1} N (d, 1,1,1 > e_{2,2} 
) = d^2 e_{2,1} N_{e_{2,2}}$ (the latter equality due to the fact 
that determining the position of $p_1, q_1, p_2$ determines that of 
$q_2$ via cross-ratio).  Thus one is finally reduced to computing, e.g.,
$N(d, 1, 1 > (1))$.  Let's identify $(p_1, q_1) = (0, \infty )$, which
then identifies cross-ratio with ordinary ratio.  A moment's reflection
shows that $N(d, 1, 1 > (1)) = d^2 M_d$ where $M_d$ is the number of 
pairs $(a,b) \in \C \times \C$ such that $a/b \in \{ \lambda_1, \lambda^{-1}\}$
for a fixed general $\lambda \in \C$ and for a general degree-$d$ map 
$f: \P^1 \to \P^1, \ \ f(a) = f(b)$.  Specializing $\lambda \to 1$ keeping 
track of multiplicities and using the Riemann-Hurwitz formula, it is 
elementary that $M_d = 4 (d-1)$, so   
$$
N(d, 1, 1 > (1)) = 4(d-1) d^2 . 
\tag2.d.5
$$
Similarly,
$$
N(d, (1) > (1)) = 2 (d-1) d^2 . 
\tag2.d.6
$$
(here one divides by 2 due to the involution $(p_1, q_1)(p_2 q_2)$).  Thus 
one can compute all the above-mentioned cross-ratio counts. 
\ss
\subheading{\bf 3. Higher genus}

\ss 
For plane curves of positive genus, there are (at least) two types of 
counts one may wish to carry out, depending on whether the moduli of the 
curve are fixed or unrestricted.  For unrestricted moduli on has the number 
$N(d,g)$ of (integral) curves of degree $d$ and geometric genus $g$ 
through $3d+ g - 1$ general points, for which a recursive formula was 
given in our earlier paper [R].  For fixed moduli, one has the 
number $N(d, g\rangle$ of integral curves of degree $d$ birational to a fixed 
general smooth curve of genus $g$ and passing through $3d-2g+2$ general
points if $g \geq 2$ (or $3d-1$ points if $g=1$).  The case $g=1$ was 
done by Pandharipande [P], who shows
$$
N(d, 1 ) = \frac{(d-1) (d-e)}{2} N(d,0) . 
\tag3.1
$$
We now show how the method of this paper yields easily a procedure for
computing $N(d_1, 2)$; see [R$'$]  
for a different approach to $N(d,g\rangle$ in general.

Now specializing an abstract curve $C$ of genus 2 a general binodal
rational curve $C_0$, with normalization $(\P^1, \{p_1,q_1\} , 
\{p_2, q_2 \}) \to (C_0, \text{node}, \text{node})$, we see that 
$N(d,2\rangle$ may be identified with the number of maps
$$
f : \P^1 \to \P^2
$$
with image $\bar{C} = f(\P^1)$ passing through $3d-2$ general points,
such that, for a fixed quadruple $(p_1, q_1, p_2, q_2), f(p_i) = 
f(q_i), i = 1,2$, up to identifying $f \sim f \o \alpha$ where $\alpha$ is 
the unique projective automorphism inducing the permutation 
$(p_1, q_1) (p_2, q_2)$ (i.e. the map induced by the limit of the 
hyperelliptic involution on $C$).  To this Theorem 1 is applicable, 
and if remains to list the boundary components $Z$.  These come in two 
types:
\def\cupp{\operatornamewithlimits{\cup}}

(1)  $Z_{d_1, d_2}^1$, which corresponds to maps
$$
(\P_1^1, \{p_1^0, q_1\} , \{ p_2, q_2\} ) \cupp\limits_{p_1^0} 
(\P_2^1 , p_1, p_1^0) 
{\overset (f_1, f_2) \to \rightarrow} \P^2
$$
where $im \ f_i = \bar{C}_i$ has degree $d_i, f_1 (p_1^0) = f_2 (p_1^0),
f_1 (q_1) = f_2 (p_1), f_1 (p_2) = f_1(q_2)$.  For enumerate
$Z_{d_1, d_2}^1$ we introduce nodal crop-ratio counts 
$$
N(d < (e)) = \# \{f: (\P^1, \{p_1, q_1\} , \{p_2, q_2\} ) \to (\P^2, 
E, \text{point}) \}
$$
where $E$ is a fixed curve of degree $e$, the point is unspecified (i.e.,
the conditions is $f(p_2) = f(q_2))$ and $im \ (f)$ is a degree-$d$ curve 
through $3d-2$ general points; similarly $N(d > (e))$ where $im\ f$ 
(and the cross-ratio) are fixed and $E$ is rational degree-$e$ through 
$3e-2$ general pints.  With these we have, e.g.,
$$
d^{1,1} (Z_{d_1,d_2}^1) = \pmatrix 3d-2 \\ 3d_2 -2\endpmatrix N_{d_2}
N(d_1 < (d_2)) + \pmatrix 3d-2\\ 3d_1 -2 \endpmatrix N_{d_1}
N(d_1 > (d_2)) . 
\tag3.2
$$
The numbers $N(d < (e))$ and $N(d > (e))$ may be computed recursively in 
analogy --as well as linkage with the cross-ratio numbers of \S2.d:
e.g. the recursion for $N(d < (e))$ involves $N(d_1 < (e))$ when 
$p_2, q_2$ go to a node on one component, as well as $N(d_1 < (d_2), (e))$ 
where the $\P^1$ splits so $p_1$ goes off to another component; and 
then will also be a $N(d_1, (e) > (d_2))$ where  an unmarked curve of 
degree $d_2$ varies.  Details are similar to the above.
$$
Z_{d_1, d_2}^2 = \left\{(\P_0^1, \{p_1^0, q_1^0\}, \{p_2, q_2\}) 
\cupp\limits_{p_1^0}
(\P_1^1, p_1, p_1^0) \cupp\limits_{q_1} (\P_2^1, q_1, q_1^0 )
{\overset (f_0, f_1, f_2)\to \rightarrow } \P^2
\right.
\tag2
$$
where $\bar{C}_i = f_i (\P^1_i)$ have degrees $d_0 = 0,d_1, d_2, \bar{C}_1$ 
and $\bar{C}_2$ are {\it tangent} at $\bar{C}_0
= f_1 (p_1^0) = f_2 (q_1^0)$ and meet at $f_1 (p_1) = f_2 (q_1)$.  This 
component has $\ell = 2$, and may be easily enumerated as in \S2, e.g.
$$
\align
d^{1,1}(Z_{d_1, d_2}^2) &= B_{d_1, d_2, 0} \pmatrix 3d-2\\ 3d_1-2\endpmatrix
. (d_1 d_2 - 2) + B_{d_2, d_1, 0} \pmatrix 3d-2 \\ 3d_1 -2\endpmatrix
(d_1 d_2 -2 ) . \\
&= (2 (d_2 -1) N_{d_1} + d_2 B_{d_1} ) \pmatrix 3d-2\\ 3d_2 - 2\endpmatrix
(d_1 d_2 - 2) \\
&+ (2(d_1 -1) N_{d_2} + d_1 B_{d_2})
\pmatrix 3d-2\\ 3d_1 - 2 \endpmatrix (d_1 d_2 -2 ) . 
\tag3.3
\endalign
$$
In this way we may obtain recursions for all the nodal cross-ratio numbers
and hence compute $N(d, 2 \rangle$.

\ss
\heading{\bf Appendix:  Rational curves on Del Pezzo surfaces}
\endheading
\ss
The following result, though quite natural, seems to our knowledge to have 
escaped explicit mention in the literature.

\proclaim{Proposition A}  Let $S$ be a Del Pezzo surface of anticanonical 
degree at least 3 and $c$ a divisor
class on $S$ with $c^2 > 0 > c.K$.  Then $c$ contains an irreducible rational
curve.
\endproclaim

\demo{Proof}  We consider the case of a blown-up plane.  By Riemann-Roch,  
$c$ is clearly effective.  If the (anticanonical) degree $-K.c=1, 2,$
our claim is easy, and likewise if $c = - K$.  In general, we use an 
induction on $-K.c$.  Let $D$ be a line on $S$ (i.e., $D^2 = - 1 = D.K$)
such that $c.D$ is minimal, and set $c' = c -D$.  If $c.D < 0$, 
induction clearly applies to $c$'s yielding an irreducible rational 
representative $C'$ moving in a family of rational curves of dimension
$-K.c-2\geq 1$, which $C'$ can therefore be specialized to (an image of)
a connected rational chain $C_0$ meeting $D$.  Then by easy and standard deformation theory (using ampleness of $-K$), $C_0 + D \simeq c$ can be 
deformed to an irreducible rational curve.  If $c.D = 0$, look instead for 
the smallest strictly positive $c.D', D' = \text{line}$ and argue as below.
(Or alternatively blow $D$ down and restart the argument-which basically
amounts to the same thing.)
Now let's assume $c.D > 0$.  It will suffice to prove 
$$
(c-D)^2 > 0 
$$
for then deformation theory can be applied as above to conclude.  To 
this end let us consider a suitable model of $S$ as a blowup of the plane
in $r$ points, $r \leq 6$, and represent $c$ in the usual way as  
$$
(b; a_1 \geq \cdots \geq a_r)
$$
where we may assume 
$$
a_r = c.E_r = c.D
$$
is minimal as above, which translates into
$$
b \geq a_1 + a_2 + a_r
$$
$$
2b \geq a_1 + \cdots + a_5 + a_r, 
$$
and our claim amounts to 
$$
N (c) := b^2 - a_1^2 - \cdots - a_r^2 > 2 a_r + 1. 
$$
Let's do some tail-flattening on the sequence 
$a_1, \cdots a_r $.  Now if $a_r < a_{r-1} < a_2 $, say, we may 
perform a `switching transformation' with
$$
a_{r-1}' = a_{r-1} - 1 , a_{r-2}' = a_{r-2} + 1
$$
which will only decrease $N(c)$.  Continuing in this way we may 
eventually achieve
$$
a_3 = \cdots = a_r
$$
Then doing a similar `switch' with $a_2$ and $a_1$ we may further 
assume $a_2 = a_r$.  Now suppose $a_1 \leq 2a_r - 1$.  Clearly
$a_r \leq 1/3 b$, so
$$
b^2 - a_1^2 - \cdots a_r^2 \geq \frac{6-r}{9} 6^2 + 4 a_r - 1 
$$
which is $> 2a_r + 1$ unless $a_r = 1$ in which case our assertion 
clearly holds anyway except for the canonical class $c=-K = (3;1, 
\cdots 1)$ (and hence would hold for our `original' class $c$
if different from $-K$).  Now if on the other hand $a_1 \geq 2a_r$, we 
have
$$
b^2 \geq (a_1 + 2a_r)^2 \geq a_1^2 + 12 a_r^2 > a_1^2 + 
(r-1) a_r^2 + 2 a_r + 1
$$
as $r\leq 6$. \qed
\enddemo

\ls
\heading{\bf References}\endheading

\ss
\roster
\item"{[FP]}"  W. Fulton, R. Pandharipande:  `Notes on stable maps and 
quantum cohomology', to appear.
\ss
\item"{[H]}"  J. Harris:  `On the Severi problem' Invent. Math. {\bf 84}
(1986), 445--461.
\ss
\item"{[P]}"  R. Pandharipande:  `A note on elliptic plane curves with 
fixed $j$-invariant' (preprint).
\ss
\item"{[R]}"  Z. Ran:  `Enumerative geometry of singular plane curves'.
Invent. Math. {\bf 97} (1989), 447-469.
\ss
\item"{[R$'$]}"  $\underline{\hskip.25in}$:  (in preparation).

\endroster

\enddocument